\documentclass[11pt]{amsart}
\usepackage{amssymb,amsmath,amsfonts,amssymb, latexsym, verbatim, esint,amsthm, mathrsfs}
\newcommand{\be}{\begin{equation}}
\newcommand{\ee}{\end{equation}}
\newcommand{\la}{\label}
\newcommand{\ba}{\begin{array}{l}}
\newcommand{\ea}{\end{array}}
\newcommand{\Rr}{{\mathbb R}}

\newcommand{\pa}{\partial}
\newcommand{\fr}{\frac}
\newcommand{\na}{\nabla}
\newtheorem{thm}{Theorem}
\newtheorem{prop}{Proposition}
\newtheorem{lemma}{Lemma}
\newtheorem{defi}{Definition}

\newcommand{\D}{\mathcal D}
\newcommand{\h}{L^2(\mathbb R^2)}
\newcommand{\bl}{Lim_{t\to\infty}}
\title[Absence of anomalous dissipation]{Absence of anomalous dissipation of energy in forced  two dimensional fluid equations}
\author{Peter Constantin}
\address{Department of Mathematics, Princeton University, Princeton, NJ 08544}
\email{const@math.princeton.edu}

\author{Andrei Tarfulea}
\address{Department of Mathematics, Princeton University, Princeton, NJ 08544}
\email{tarfulea@math.princeton.edu}

\author{Vlad Vicol}
\address{Department of Mathematics, Princeton University, Princeton, NJ 08544}
\email{vvicol@math.princeton.edu}
\date{}
\begin{document}
\maketitle
\tableofcontents
\section{Introduction}
Anomalous dissipation of energy in three dimensional turbulence is one of the basic statements of physical theory ~\cite{Frisch95}. It has been verified experimentally to a large degree ~\cite{SreenivasanAntonia97}, but not mathematically. The statement is about the average behavior of the energy dissipation rate
\[
\epsilon = \nu\langle |\na u|^2\rangle
\]
as $\nu\to 0$. Here $\nu$ is kinematic viscosity, $u$ is the velocity (assumed to have mean zero), $\na$ are spatial gradients and $\langle\dots\rangle$ represents an ensemble average or space-time average. The assertion of turbulence theory is that $\epsilon$  is a {\em{positive number}}, and that it does not vanish with viscosity, in the limit of zero viscosity. The term ``anomalous dissipation'' was imported from field theory in physics and it refers to the fact that, in the limit of vanishing viscosity, there still is remanent dissipation, even though the limit equation conserves energy. 

There are two distinct approaches to the question of anomalous dissipation. In the first, the limit of zero viscosity is taken on solutions of the initial value problem with fixed initial data. Under appropriate conditions this leads to a solution of the corresponding initial value problem of the inviscid equation. This equation conserves energy if solutions are smooth, but might dissipate energy if solutions are not sufficiently smooth. This circle of ideas, and specifically the precise degree of smoothness needed, goes by the name of ``Onsager conjecture'' ~\cite{EyinkSreenivasan06,Eyink94,Eyink05,ConstantinETiti94,CheskidovConstantinFriedlanderSvydkoy08,DeLellisSzekelihidi09,DeLellisSzekelihidi10,DeLellisSzekelyhidi12,DeLellisSzekelihidi12b,DeLellisSzekelihidi12c,Isett12}. This approach is therefore about the initial value problem for the limit equations and it requires lack of smoothness of solutions. The blow up problem is open for 3D incompressible Euler equations, and this allows to envision the possibility of existence of dissipative solutions arising from smooth initial data. Anomalous dissipation of energy  can be proven for incompressible 2D Euler equations as well, for very rough solutions, although in 2D non-smooth solutions cannot arise spontaneously from smooth ones. The class of dissipative solutions of the inviscid equations is very large indeed.

The second way of looking at the anomalous dissipation issue is to take long time averages first, in order to achieve a ``permanent regime'' of the viscous equations, and only then send the viscosity to zero.  This second approach is espoused in this paper. Denoting by $S^{(\nu)}(t,u_0)$ the solution of the viscous equation at time $t\ge 0$ which started at time $t=0$ from the initial data $u_0$, the second approach looks at
\[
\langle |\na u|^2\rangle = \lim_{T\to\infty}{\fr{1}{T}}\int_0^T\int_{\Rr^d}\left|\na(S^{(\nu)}(t,u_0))\right |^2dxdt
\]   
and asks if $\lim_{\nu\to 0}\nu\langle|\na u|^2\rangle = \epsilon$ is positive or not.  A significant question is that of dependence of forcing and initial data. In the absence of forcing the long time limit vanishes even for Leray weak solutions of 3D Navier-Stokes equations: the ``permanent regime'' is trivial, and turbulence is decaying. One must then take finite time averages, with times of durations that diverge with vanishing viscosity, but not too fast, nor too slow. This unforced case is perhaps the only case in which a general global, a priori upper bound on $\epsilon$ that is viscosity independent is rigorously known.

The long time averaged dissipation has a chance of being not trivial if the flow is forced, either at boundaries or if body forces stir the flow. A conceptual difficulty arises then because there exist situations in which the rate of dissipation, as defined, is infinite.  There are rigorous studies ~\cite{DoeringFoias02,DascaliucFoiasJolly09,DascaliucFoiasJolly10} where bounds for $\epsilon$ are expressed in terms of the average kinetic energy of the solutions in the case of forced Navier-Stokes equations; however, there are no viscosity  (or Reynolds number) -independent a priori bounds on the average kinetic energy. The question of 
obtaining examples and a characterization of flows with uniform upper bounds for $\epsilon$ is open. In fact, the opposite situation can be easily found: $\epsilon$ can be unbounded when we consider spatially periodic 2D forces that are eigenfunctions of the Stokes operator. We write the Navier-Stokes equation symbolically as
\[
\pa_t u + \nu Au + B(u,u) =f
\]
where $A= -{\mathbb{P}}\Delta$ is the Stokes operator with $\mathbb P$ the projector on divergence-free vectors and $B(u,u) = {\mathbb P}(u\cdot\na u)$ is
the quadratic nonlinearity. We take   $Af=\lambda f$ with $f$ and $u$  periodic, divergence-free functions of two spatial variables. We consider  $u_f = \fr{1}{\nu\lambda }f$. This is a smooth, time independent solution. Indeed, $u_f$ satisfies the time-independent, unforced Euler equation $B(u_f, u_f)=0$, and, of course $\nu Au_f =f$, as well.  If $u_0=u_f$
then $\langle |\na u_f|^2\rangle = \nu^{-2}\lambda^{-1}\|f\|_{L^2}^2$, where 
$f$ does not depend on $\nu$. This makes the limit of $\epsilon$ infinite as viscosity vanishes. The steady solutions $u_f$ are perfectly admissible as 3D periodic flows as well. They are unstable if $\lambda$ is not the first eigenvalue of the Stokes operator, but stable in 2D if $\lambda$  is the first eigenvalue. Incidentally, the upper bound $\epsilon\le C(U^3\sqrt{\lambda} + \nu \lambda U)$ (with $U =\|u_f\|_{L^2}$) of ~\cite{DoeringFoias02} is true in this case as well, and it is imprecise, with both left and right hand side diverging as $\nu\to 0$, but at different rates.
An interesting recent asymptotics and numerical study  ~\cite{GalletYoung13} reports finding solutions of the 2D Navier-Stokes equation that ``settle'' to a condensate that has a nontrivial component in the first eigendirection of the corresponding Stokes operator and has bounded amplitude as viscosity vanishes. This of course is impossible for all initial data, as demonstrated above, but it is an intriguing possibility for parts of the phase space. It is known that if we assume  that an initial data $u_0$ is smooth enough then the solution of the Navier-Stokes equations with smooth forcing (even if not an eigenfunction of the Stokes operator) converges to the corresponding solution of the forced Euler equation $u(t) = S^{(0)}(t,u_0)$ on $[0,T]$, for any $T$, a function that solves 
\[
\pa_t u + B(u,u) =f
\]
with initial data $u_0$. For Kolmogorov forcing, (forces which are eigenfunctions of the Stokes operator), the putative existence of time independent solutions $u^{(\nu)}$ which are  uniformly  bounded in $\nu$ in energy norm, implies the convergence of (a subsequence of) $u^{(\nu)}$ to a time independent solution $u^{0}$ of the forced incompressible Euler equations, $B(u^0,u^0) =f$. If the solutions  $S^{(\nu)}(t,u_0)$ are at a bounded distance from $u^{(\nu)}$  uniformly in time, for large time, one can prove  that $S^{(0)}(t,u_0)$ are at a the same bounded distance from $u^0$ for large time. In particular, if $S^{(\nu)}(t,u_0)$ converge uniformly in time to $u^{(\nu)}$, then $S^{(0)}(t,u_0)$ converge in time to $u^{(0)}$. Smooth steady solutions $B(u^{(0)}, u^{(0)}) =f$ of forced Euler equations with Kolmogorov forcing can be easily constructed, but for which initial data solutions converge to them is another matter. Such behavior, if it exists at all, must be rather special.

The dynamics of the forced Euler equation, the existence of bounded sequences of stationary solutions of the periodic, forced Navier-Stokes equations, and even of solutions with bounded average dissipation of energy are open problems.

Bounds can be obtained for 2D forced Navier-Stokes equations with  bottom drag
(friction),
\[
\pa_t u + \gamma u + \nu Au + B(u,u) =f
\]
where the friction coefficient $\gamma>0$ is kept fixed. Then the energy is bounded (in terms of $\gamma$) uniformly for small $\nu$ and actually even the
enstrophy ($H^1$ norm) is bounded uniformly in $\nu$. Consequently, there is no anomalous dissipation of energy. The absence of anomalous dissipation of enstrophy is more subtle because there are no upper bounds for the average $H^2$ norm, for arbitrary forces. The paper ~\cite{ConstantinRamos07} proves nevertheless that the dissipation of enstrophy vanishes in the limit of zero viscosity, for arbitrary time-independent forces. 

In this paper we prove absence of anomalous dissipation of energy for 
surface quasi-geostrophic (SQG) equations. These equations have generated a lot of attention in recent years ~\cite{ConstantinMajdaTabak94,HeldPierrehumbertGarnerSwanson95,Resnick95,OhkitaniYamada97,Cordoba98,ConstantinNieSchorghofer98,ConstantinWu99,CordobaFefferman01,ConstantinCordobaWu01,Wu01,Constantin02,Berselli02,CordobaFefferman02,ChaeLee03,SchonbekSchonbek03,CordobaCordoba04,Ju04,Wu04,Ju05a,CordobaFontelosManchoRodrigo05,DengHouLiYu06,Miura06,Ju07,ChenMiaoZhang07,HmidiKeraani07,KiselevNazarovVolberg07,AbidiHmidi08,Chae08,ConstantinWu08,DongDu08,Marchand08a,DongLi08,Yu08,ConstantinWu09,DongPavlovic09,DongPavlovic09b,KiselevNazarov09,FriedlanderPavlovicVicol09,CaffarelliVasseur10,CaffarelliVasseur10b,KiselevNazarov10,Dong10,OhkitaniSakajo10,Silvestre10a,Dabkowski11,FeffermanRodrigo11,Kiselev11,ConstantinVicol12,FeffermanRodrigo12,DabkowskiKiselevVicol12,XueZheng12,DabkowskiKiselevSilvestreVicol12,ConstantinLaiSharmaTsengWu12}. \\ We are interested in the question of anomalous dissipation for forced, viscous critical SQG. We consider the equation
\[
\pa_t \theta + (R^{\perp}\theta)\cdot\na\theta + \gamma\D\theta -\nu\Delta \theta =f
\]
in $\Rr^2$, where $\D = {\mathbb I} + (-\Delta)^{\fr{1}{2}}$ is the damping operator, $R^{\perp} = (-R_2,R_1)$ are Riesz transforms, $f\in L^{\infty}(\Rr^2)\cap L^1(\Rr^2)$ is time independent deterministic forcing, $\gamma>0$ is fixed and $\nu>0$. We prove that there is no dissipative anomaly,
\[
\lim_{\nu\to 0}\nu\langle |\na\theta |^2\rangle = 0
\]
where $\langle\cdots\rangle$ is space-time average on solutions. The proof of absence of anomalous dissipation follows the same blueprint as the proof in ~\cite{ConstantinRamos07}. We establish first that the viscous  semi-orbits are relatively compact in the phase space. Then we introduce the adequate statistical solutions for both viscous ($\nu>0$) and inviscid ($\nu =0$) equations. These are measures in phase space, arising naturally as long time limits on solutions. The next step is to prove that the zero viscosity limits of statistical solutions of the viscous equations are statistical solutions of the inviscid equations, and that these preserve the energy balance. Once this is achieved, the absence of anomalous dissipation follows by an argument by contradiction. There are a number of technical difficulties encountered in the proof for SQG that are not present in the case of 2D Navier-Stokes. In order to obtain the uniform integrability property on positive semiorbits we use nonlocal calculus identities. The weak continuity of the nonlinear term is proved using a commutator structure of the nonlinearity, a structure that was used already in ~\cite{Resnick95}. The energy balance is proved using a formula for nonlinear fluxes ~\cite{ConstantinETiti94} and a bound in $H^{\fr{1}{2}}$ that is available for critical SQG, and that replaces the Besov space argument of ~\cite{ConstantinETiti94,CheskidovConstantinFriedlanderSvydkoy08}.

The rest of the paper is organized as follows. In section 2 we make more precise the comments about Kolmogorov forced Navier-Stokes and Euler equations. In section 3 we present the forced viscous SQG equations and prove some properties of solutions, including the relative compactness of positive semiorbits. In section 4 we introduce the notion of stationary statistical solutions of the viscous equations. In section 5 we prove that inviscid limits of stationary statistical solutions are stationary statistical solutions of the forced
critical SQG equations which preserve the energy dissipation balance. In section 6 we construct stationary statistical solutions using time averages and in section 7 we present the argument by contradiction and concluding remarks.
\section{2D Forced Navier-Stokes Equations}
We consider 2D periodic incompressible Navier-Stokes equations
\[
\pa_t u -\nu\Delta u + u\cdot\na u +\na p = f
\] 
where $u:[-\pi L,\pi L]^2\times [0,\infty)\to \Rr^2$ is divergence-free, $\na\cdot u =0$, and periodic, $u(x\pm 2\pi Le_i,t) = u(x,t)$ (here $e_i$, $i=1,2$ is the canonical basis of $\Rr^2$). We take time independent $f:[-\pi L,\pi L]^2\to \Rr^2$ that is divergence-free $\na\cdot f=0$, periodic of the same period $2\pi L$, $f(x\pm 2\pi Le_i) =f(x)$, and an eigenfunction of the Stokes operator, which in the case of divergence-free periodic function is just the Laplacian on each component, $-\Delta f = \lambda f$. We refer to such forcing as ``Kolmogorov forcing''.
We choose to measure lengths in units of $L$, and because the force plays an important role and has units of $[f]={\rm{length}}\times{\rm{time}}^{-2}$,  we measure time in units of $T=\sqrt{\fr{L}{F}}$ where $F$ is the RMS force, $F^2=(2\pi L)^{-2}\int_{|x_i|\le \pi L,\, i=1,2}|f(x)|^2dx$. Rescaling, i.e considering $u = \fr{L}{T}{\widetilde{u}}(\fr{x}{L}, \fr{t}{T})$, $f = F{\widetilde{f}}(\fr{x}{L})$, $p = \fr{L^2}{T^2}{\widetilde{p}}(\fr{x}{L}, \fr{t}{T})$ and $\nu = \fr{L^2}{T}\widetilde{\nu}$, and dropping tildes, we have thus
\be
\pa_t u -\nu \Delta u + u\cdot\na u + \na p = f, \quad \na\cdot u = 0
\la{nse}
\ee
with $u:[-\pi,\pi]^2\times [0,\infty)\to \Rr^2$, $f:[-\pi,\pi]^2\to \Rr^2$ of period $2\pi$, with normalized $L^2$ norm equal to 1, and $\nu$ nondimensional, in fact the inverse Reynolds number. We still have $\na\cdot f = 0$ and
\be
-\Delta f = \lambda f\la{kol}
\ee
with the nondimensional (new) $\lambda$ equal to the dimensional (old) $\lambda $ multiplied  by $L^2$. The Fourier series representation of $u$ is
\be
u(x,t) = \sum_{j\in {\mathbb Z}^2}\widehat{u}(j,t) e^{ij\cdot x}
\la{uf}
\ee
with ${\widehat{u}}: {\mathbb Z}^2\times [0,\infty)\to {\mathbb C}^2$. Without loss of generality the average of $u$ vanishes, ${\widehat{u}}(0,t) =0$. Because $\na\cdot u =0$ and we are in two dimensions, without loss of generality
\be
{\widehat{u}}(j,t) = u_j(t) \fr{j^{\perp}}{|j|}\la {uj}
\ee
where $j^{\perp} = (-j_2,j_1)^*$ and $v^*$ is the transpose. Now $u_j(t)$ is a scalar complex valued function of time, and the requirement that $u$ be real valued implies the requirement that ${\overline{u_j}} = - u_{-j}$ (from $\overline{\widehat{u}(j)} = {\widehat{u}}(-j)$).
We note that the stream function, defined by the relation $u = \na^{\perp}\psi$
\[
\psi (x,t) = \sum_{j\in {\mathbb Z}^2}{\widehat{\psi}}(j,t) e^{ij\cdot x},
\]
has Fourier coefficients
\[
{\widehat{\psi}}(j,t) = -i |j|^{-1} u_j,
\] 
or, in other words $u_j = i|j|\widehat{\psi}(j)$. If $u$ is divergence-free, it does not necessarily follow that $u\cdot\na u $ is divergence-free as well. The projector on divergence-free functions is computed for 2D Fourier series 
\[
v(x) = \sum_{j\in {\mathbb{Z}}^2\setminus\{0\}}{\widehat{v}}(j)e^{ij\cdot x}
\]
as
\[
{\mathbb Pv}(x) =  \sum_{j\in {\mathbb{Z}}^2\setminus\{0\}}{\mathbb P}_j{\widehat{v}}(j)e^{ij\cdot x}
\]
with
\[
{\mathbb P}_j v = \left(v\cdot\fr{j^{\perp}}{|j|}\right)\fr{j^{\perp}}{|j|}.
\]
The Stokes operator, denoted $A$, is 
\be
A = -{\mathbb P}\Delta
\la{stok}
\ee
and the projection of the bilinear term is
\be
B(u,v) = {\mathbb P}(u\cdot\na v).
\la{buv}
\ee
Using our convention that mean-free, divergence-free vectors are written as
\[
v(x) = \sum_{j\in{\mathbb Z}^2\setminus\{0\}}v_j\fr{j^{\perp}}{|j|}e^{ij\cdot x}
\]
with $v_j$ complex scalars, we obtain for divergence-free $u$ and $v$,
\be
\left[B(u,v)\right]_l = i\sum_{j+k=l, \; j,k,l \neq 0}u_jv_k\left (\fr{j^{\perp}}{|j|}\cdot k\right)\left(\fr{k^{\perp}}{|k|}\cdot \fr{l^{\perp}}{|l|}\right).
\la{buvl}
\ee
In particular
\be
\left[B(u,u)\right]_l = i\sum_{j+k=l, \; j,k,l \neq 0}u_ju_k\left(j^{\perp}\cdot k\right)\left (k\cdot l\right )\fr{1}{|j||k||l|}
\la{buul}
\ee
and because of the antisymmetry of $u_ju_k\left(j^{\perp}\cdot k\right)\fr{1}{|j||k|}$ in $j,k$ at fixed $l$, we have
\be
\left[B(u,u)\right]_l = \fr{i}{2}\sum_{j+k=l, \; j,k,l \neq 0}u_ju_k\left(j^{\perp}\cdot k\right)\left (|k|^2-|j|^2\right )\fr{1}{|j||k||l|}
\la{buuls}
\ee
This shows that the only contributions to $B(u,u)$ come from distinct energy shells , i.e. $|j|\neq |k|$. In particular, any function whose Fourier support is on a single energy shell, solves $B(u,u) =0$. This is the case for eigenfunctions of the Stokes operator.  In terms of the vorticity, if $u = \nabla^{\perp}\psi$ and $\Delta \psi = \lambda \psi$, it follows that $u\cdot\nabla\omega =0$
because the vorticity $\omega = \nabla^{\perp}\cdot u$ is given by $\omega= \Delta \psi$. The 2D incompressible unforced Euler equation can be written in vorticity formulation as
\[
\pa_t \omega + u\cdot\na \omega = 0
\]
and therefore, if $\Delta \psi = \lambda \psi$, we obtain time independent 
solutions of the Euler equations. Another way of seeing that eigenfunctions of the Stokes operator are steady solutions of unforced Euler equations  is via the identity
\be
AB(u,u) = B(u,Au)-B(Au,u).\la{idb}
\ee
This is proven by observing that for 2D divergence-free vectors $u$, 
\[
\Delta(u\cdot\nabla u_i)- u\cdot\nabla \Delta u_i + \Delta u\cdot\nabla u_i= 2\pa_k(({ det}\na u) \delta_{ik}).
\]
Thus, if $Au =\lambda u$ then $B(u,u) =0$, because $A$ is invertible. In particular, if $Af = \lambda f$ then the time-independent $u=u_f$ given by  $u_f = \fr{1}{\nu\lambda }f$  solves the Navier-Stokes equation (\ref{nse}). Let us consider now solutions $u(t)$ of the initial value problem (\ref{nse}) with divergence-free smooth initial data (it is enough to consider $H^1$ initial data). These are unique, exist for all time, become instantly infinitely smooth, and converge in time to a compact, finite dimensional attractor ~\cite{ConstantinFoias88}. The attractor contains $u_f$ and its unstable manifold. In particular, it follows that the largest norm of functions in the attractor (any norm) diverges with $\nu$. If the diameter of the attractor would be bounded, then $\epsilon$ would diverge as $\nu\to 0$, for any space time average on trajectories. 

Let us remark that if we fix smooth, divergence-free initial data $u_0\in H^s$, $s>2$ then  
\[
\lim_{\nu\to 0}S^{(\nu)}(t,u_0) = S^E(t,u_0)
\]
holds where $S^{(0)}(t,u_0)$ is the unique global solution of 
\be
\pa_t u + B(u,u) =f
\la{feq}
\ee
with initial data $u_0$. The convergence is in $C([0,T], H^{s'})$,  $s'<s$, for any $T$.
This follows from the global existence of smooth solutions of the forced Euler equations and from convergence as long as these solutions are smooth ~\cite{BealeKatoMajda84,Constantin86,Masmoudi07}.
This result does not need $f$ to be an eigenfunction  of the Stokes operator, only to be smooth enough. The long time behavior of $S^{(0)}(t,u_0)$ and that of $S^{(\nu)}(t,u_0)$ for small $\nu$ can be very different. In fact, if $f=0$, the behavior {\em{is}} different, because the inviscid solution conserves the initial energy, while the viscous solution converges to zero. 

Let us consider now Kolmogorov forcing and any family of steady solutions
$u^{(\nu)}$ of the forced Navier-Stokes equations
\be
\nu Au^{(\nu)} + B(u^{(\nu)},u^{(\nu)}) = f.
\la{stns}
\ee
Taking the scalar product with $u^{(\nu)}$ and then with $Au^{(\nu)}$, we have
\[
\nu \|u^{(\nu)}\|_{H^1}^2 = (f,u^{(\nu)})_{L^2}
\]
and
\[
\nu \| Au^{(\nu)}\|^2_{L^2} = \lambda (f,u^{(\nu)})_{L^2}
\]
where we used the notation
$(f,u)_{L^2} = \fr{1}{(2\pi)^2}\int_{[-\pi,\pi]^2}(f\cdot u) dx$ and the facts that $\|u\|^2_{H^1} = (u, Au)_{L^2}$,  $Af= \lambda f$, $(B(u,u), u)_{L^2}= 0$, and  $(B(u,u), Au)_{L^2}=0$. Subtracting we have
\be
\|Au^{(\nu)}\|^2_{L^2} = \lambda \|u^{(\nu)}\|^2_{L^2}.
\la{h2bal}
\ee
Using the  straightforward inequality
\[
\fr{\|u\|_{H^1}^2}{\|u\|_{L^2}^2} \le \fr{\|Au\|_{L^2}^2}{\|u\|_{H^1}^2}
\]
and {\em{assumimg}} that the family is uniformly bounded in $L^2$:
\be
\|u^{(\nu)}\|^2_{L^2}\le E,
\la{steb}
\ee
it follows that
\be
\|Au^{(\nu)}\|_{L^2}^2 = \lambda \|u^{(\nu)}\|_{H^1}^2 \le \lambda^2 E.
\la{stbounds}
\ee
Now we can pass to a convergent subsequence, first weakly convergent in $L^2$, but because of compact embedding of $H^1$, strongly in $L^2$, and by the same argument, weakly in $H^2$ and strongly in $H^1$. There is therefore enough control to show that the limit $u^{(0)}$ is a steady solution of the forced Euler equations,
\be
B(u^{(0)}, u^{(0)}) = f.
\la{steu}
\ee
Similarly, for time dependent solutions of (\ref{nse}), $u(t) = S^{(\nu)}(t,u_0)$, we bound the difference $\|u\|_{H^1}^2 -\lambda \|u\|_{L^2}^2$.  Indeed, the evolution of the $L^2$ and $H^1$ norms are given by
\be
\fr{d}{2dt}\|u\|_{L^2}^2 + \nu \|u\|_{H^1}^2 = (f,u)_{L^2}
\la{l2ev}
\ee
and
\be
\fr{d}{2dt}\|u\|_{H^1}^2 + \nu \|Au\|_{L^2}^2 = \lambda(f,u)_{L^2}
\la{h1ev}
\ee
and subtracting we have
\be
\fr{d}{2dt}\left[ \|u\|_{H^1}^2-\lambda\|u\|_{L^2}^2\right] + \nu\left [
\| Au\|_{L^2}^2 -\lambda \|u\|^2_{H^1}\right] = 0.
\la{diffeq}
\ee
Let us denote 
\be
\delta(t) =  \|u\|_{H^1}^2-\lambda\|u\|_{L^2}^2
\la{delta}
\ee
and
\be
\mu(t) = \fr{\|u\|_{H^1}^2}{\|u\|_{L^2}^2}.
\la{mu}
\ee
Let us observe that
\[
\mu(t)\ge \lambda_1 =1
\]
where $\lambda_1$ is the smallest eigenvalue of $A$, and that
\[
\|Au\|^2_{L^2}-\mu^2\|u\|_{L^2}^2 = \|u\|_{L^2}^2\left \| (A-\mu)\fr{u}{\|u\|_{L^2}}\right\|_{L^2}^2
\]
Adding and subtracting $\nu \|u\|_{L^2}^2\mu^2$, (\ref{diffeq}) becomes
\be
\fr{d}{2dt}\delta(t) + \nu\|u\|_{L^2}^2\left \| (A-\mu)\fr{u}{\|u\|_{L^2}}\right\|_{L^2}^2 + \nu\mu(t)\delta (t) = 0.
\la{deleq}
\ee
In particular
\be
\fr{d}{dt}\delta + 2\nu\mu\delta \le 0
\la{delineq}
\ee
and therefore  
\be
\delta(t) \le \delta(0)e^{-2\nu\int_0^t\mu(s)ds}.
\la{delt}
\ee
Note that if $\delta(0)\le 0$ then this implies that $\delta(t)\le 0$ for all $t$. If $\delta(0)>0$ then the right hand side of (\ref{delt}) decays fast to zero. In either case (\ref{delt}) shows that $\delta(t)$ is bounded on solutions, \[
\delta(t) \le \delta_{+}(0)= \max\{0,\; \delta(0)\}.
\]
This implies an automatic {\em{viscosity independent and time independent}} bound on $\|u\|_{H^1}$ {\em{given a viscosity independent and time independent bound}} on $\|u\|_{L^2}$. Let us {\em{assume}} that
\be
\sup_{\nu>0, t\ge 0 }\|S^{(\nu)}(t,u_0)\|_{L^2}^2 \le E.
\la{et}
\ee
Then, we have that
\be
\sup_{\nu>0, t\ge 0 }\|S^{(\nu)}(t,u_0)\|_{H^1}^2 \le \lambda E + \delta_{+}(0).
\la{enst}
\ee

Let us assume now that the solutions $S^{(\nu)}(t,u_0)$ have the property that
\[
\|S^{\nu}(t,u_0)- u^{(\nu)}\|_{L^2}\le \gamma
\]
for $t\ge T$, for fixed $\gamma$. Then, by passing to the limit, (on a subsequence for $u^{(\nu)}$) at each fixed $t\ge T$, we obtain that
\[
\|S^{(0)}(t,u_0)- u^{(0)}\|_{L^2}\le \gamma
\]
for all $t\ge T$. If $\gamma \to 0$ as $T\to \infty$ we obtain convergence in time of $S^{(0)}(t,u_0)$ to a solution of the steady forced Euler equations. The same thing will happen in higher norms, under the corresponding assumptions. 
It is relatively easy to construct Kolmogorov forces $f$ such that the forced, time independent Euler equation
\[
B(u,u) = f
\]
has solutions. It is enough to take two eigenfunctions $u_1$ and $u_2$ corresponding to distinct eigenvalues of the Stokes operator,
\[
Au_i = a_i u_i,\quad i=1,2,
\]
with $a_1<a_2$ and  with orthogonal spectral support, i.e. $j\perp k$ if
${\widehat u_1}(j)\neq 0$, ${\widehat u_2}(k) \neq 0$. After rotation of axes, this means $u_1$ is a function of one variable and $u_2$ a function of the orthogonal variable, e.g.
\[
u_i = \na^{\perp}\psi_i
\]
with
\[
\psi_i = \alpha_i\sin(k_ix_1) + \beta_i\cos(k_ix_i), \quad i=1,2
\]
$a_1=k_1^2<a_2= k_2^2$. Set $u = u_1+u_2$ Then $f= B(u,u)$ is an eigenfunction of the Stokes operator with eigenvalue $\lambda = a_1 + a_2$. In general $f\neq 0$.
\section{Forced, Viscous Critical SQG}
We consider the equation
\be
\pa_t \theta + u\cdot\na\theta  + \gamma{\mathcal D}\theta -\nu \Delta \theta = f
\la{nfsqg}
\ee
for a scalar valued $\theta:\Rr^2\times [0,\infty)\to \Rr$. Here 
\be
u = R^{\perp}{\theta}
\la{uth}
\ee
with $R^{\perp} = (-R_2,R_1)$, and $R= \na(-\Delta)^{-\fr{1}{2}}$ the Riesz transforms. The damping operator $\D$ is given by
\be
\D = \Lambda + 1
\la{damp}
\ee
with $\Lambda = (-\Delta)^{\fr{1}{2}}$. The coefficient $\gamma>0$ is fixed throughout the work and the coefficient $\nu>0$ is a parameter that we will let vary. The force $f\in L^1(\Rr^2)\cap L^{\infty}(\Rr^2)$ is fixed and time independent. We recall here that $\Lambda$ is defined at the Fourier transform level by
\[
{\widehat{{\Lambda \phi}}}(\xi) = |\xi|\widehat{\phi}(\xi)
\]
where
\[
\widehat{\phi}(\xi) = \int_{\Rr^2}e^{-ix\cdot\xi}\phi(x)dx
\]
and also
\[
\Lambda \phi (x) = cP.V.\int_{\Rr^2}\fr{\phi(x)-\phi(y)}{|x-y|^3}dy
\]
for an appropriate constant $c$ and smooth enough $\phi$. We will use also the pointwise identity ~\cite{CordobaCordoba04,Constantin06}
\be
2\phi(x)\cdot\Lambda\phi(x) = \Lambda (|\phi|^2)(x) + D[\phi](x)
\la{pos}
\ee
with
\be
D[\phi](x) = c\int_{\Rr^2}\fr{(\phi(x)-\phi(y))^2}{|x-y|^3}dy.
\la{dpos}
\ee 

\begin{prop}\la{nsprop} Let $\nu>0$, $f\in L^1(\Rr^2)\cap L^{\infty}(\Rr^2)$, $\theta_0\in L^{1}(\Rr^2)\cap L^{\infty}(\Rr^2)$.
The solution $\theta (x,t) =S^{(\nu)}(t,\theta_0)$ of (\ref{nfsqg}) exists for all time, is unique, satisfies the energy equation
\be
\fr{d}{2dt}\|\theta \|^2_{L^2(\Rr^2)} + \gamma\|\theta\|_{H^{\fr{1}{2}}(\Rr^2)}^2 + \nu\|\na \theta\|_{L^2(\Rr^2)}^2 =
(f,\theta)_{L^2(\Rr^2)}
\la{bal}
\ee
and the bounds
\be
\|\theta(\cdot, t)\|_{L^p(\Rr^2)} \le e^{-\gamma t}\left\{\|\theta_0\|_{L^p(\Rr^2)} -\fr{1}{\gamma}\|f\|_{L^p(\Rr^2)}\right\} + \fr{1}{\gamma}\|f\|_{L^p(\Rr^2)}\la{lpb}
\ee 
for $1\le p\le\infty$. Moreover the positive semi-orbit 
\[
O_{+}(\theta_0) =\{\theta =\theta (\cdot,t)\left |\right .\; t\ge 0\}\subset L^2(\Rr^2)
\]
is uniformly integrable: for every $\epsilon>0$, there exists $R>0$ such that
\be
\int_{|x|\ge R}|\theta(x,t)|^2dx \le \epsilon
\la{notr}
\ee
holds for all $t\ge 0$.
\end{prop}
We used the notation
\[
(f,g)_{\h} = \int_{\Rr^2}f(x)g(x)dx
\]
and we note that
\[
\|\theta \|_{H^{\fr{1}{2}}(\Rr^2)}^2 = (\D\theta, \theta)_{L^2(\Rr^2)}.
\]
The proof of existence, uniqueness and regularity follows along well established lines and will not be presented here. The bounds (\ref{lpb}) follow from the maximum principles and nonlocal calculus identities of which (\ref{pos}) is the quadratic example ~\cite{CordobaCordoba04}  and which imply that
\[
\int_{\Rr^2}\phi^{p-1}\Lambda \phi dx\ge 0
\]
if $p$ is even or if $\phi$ is nonnegative. The uniform integrability property (or ``no-travel'' property  ~\cite{ConstantinRamos07}) is proved here below. We consider the function
\[
Y_R(t) = \int_{\Rr^2}\chi\left(\fr{x}{R}\right )\theta^2(x,t)dx
\]
where $\chi$ is a nonnegative smooth function supported in $\{x\in\Rr^2\left |\right. \; |x|\ge \fr{1}{2}\}$ and identically equal to 1 for $|x|\ge 1$.
We take (\ref{nfsqg}), multiply by $2\chi(\fr{x}{R})\theta(x)$ and integrate.
The more challenging term we encounter is
\[
2\gamma\int_{\Rr^2}(\Lambda\theta(x))\chi\left(\fr{x}{R}\right)\theta(x,t)dx.
\]
Using (\ref{pos}) we have
\[
\ba
2\gamma \int_{\Rr^2}(\Lambda\theta(x))\chi\left(\fr{x}{R}\right)\theta(x,t)dx\\
\ge \gamma \int_{\Rr^2}\Lambda(\theta(x)^2)\left(1-\left(1-\chi\left(\fr{x}{R}\right)\right)\right)dx \\
=- \gamma\int_{\Rr^2}(\theta(x)^2)\Lambda\left(1-\chi\left(\fr{x}{R}\right)\right)dx
\ea
\]
where $\Lambda(1-\chi)$ is well defined because $1-\chi\in C_0^{\infty}$. Moreover $1-\chi\left(\fr{x}{R}\right) = \phi\left(\fr{x}{R}\right)$ and therefore, in view of the fact that $\Lambda(\phi(\fr{x}{R})) = \fr{1}{R}(\Lambda\phi)(\fr{x}{R})$ and $|\Lambda \phi(x)|\le C$,
\[
2\gamma \int_{\Rr^2}(\Lambda\theta(x))\chi\left(\fr{x}{R}\right)\theta(x,t)dx\ge -\fr{C\gamma}{R}\|\theta(\cdot, t)\|_{L^2(\Rr^2)}^2
\]
The contribution of the nonlinear term $u\cdot\na\theta$ is bounded by integrating by parts and using 
\[
\|u\|_{L^3(\Rr^2)}\le C \|\theta\|_{L^3(\Rr^2)}.
\] 
The contribution of the forcing term is bounded by
\[
2\left |\int_{\Rr^2} f\chi\theta dx\right |\le C\|\theta(\cdot,t)\|_{\h}\sqrt{\int_{|x|\ge \fr{R}{2}}|f(x)|^2dx}
\] 
We obtain
\[
\ba
\fr{d}{dt}Y_R(t) + 2\gamma Y_R(t)\\ 
\le \fr{C}{R}\left[\|\theta(\cdot, t)\|_{L^3(\Rr^2)}^3 + \gamma \|\theta (\cdot, t)\|^2_{\h}\right] +\fr{C\nu}{R^2}\|\theta (\cdot, t)\|_{\h}^2 \\
+ C\|\theta(\cdot,t)\|_{\h}\sqrt{\int_{|x|\ge \fr{R}{2}}|f(x)|^2dx}.
\ea
\]
Because of (\ref{lpb}) and the fact that $f^2$ is integrable, the right hand side is as small as we wish, unfiormly in time, provided $R$ is chosen large enough. The choice of $R$ depends only on $\gamma$, $f$ and on norms of $\theta_0$ in $L^2$ and $L^3$, and can be made uniformly in $\nu$ for bounded $\nu$, although we do not need this. Once we chose $R$ so that the right-hand side is less than $\gamma\epsilon$ we have the inequality
\[
Y_R(t) \le e^{-2\gamma t}Y_R(0) + \fr{\epsilon}{2}
\]
and the uniform integrability follows from the fact that $Y_R(0)$ is small for large $R$.
\section{Stationary Statistical Solutions}
We introduce first the notion of stationary statistical solution for forced viscous 
SQG, in the spirit of ~\cite{Foias72,Foias73} and~\cite{ConstantinRamos07}

\begin{defi}\la{sssn}
A stationary statistical solution of (\ref{nfsqg}) is a Borel probability measure $\mu^{(\nu)}$ on $L^2(\Rr^2)$ such that
\[
\tag{(a)} \;\; \int_{\h}\|\theta\|_{H^1}^2d\mu^{(\nu)}(\theta)<\infty
\]

\[
\tag{(b)} \; \; \int_{\h}(N^{(\nu)}(\theta), \Psi'(\theta))_{\h}d\mu^{(\nu)}(\theta) = 0
\]
for all $\Psi\in{\mathcal T}$,
and
\[
\tag{(c)} \; \; \int_{E_1\le \|\theta\|_{H^{\fr{1}{2}}}\le E_2}\left(\gamma \|\theta\|_{H^{\fr{1}{2}}}^2 + \nu\|\na\theta\|_{\h}^2 - (f,\theta)_{\h}\right )d\mu^{(\nu)}(\theta)\le 0
\]
for all $E_1\le E_2$.
\end{defi}
Here 
\be
N^{(\nu)}(\theta) = R^{\perp}\theta\cdot\na \theta + \gamma\D\theta -\nu\Delta \theta-f
\la{nntheta}
\ee
and the class of cylindrical test functions ${\mathcal T}$ is defined by
\begin{defi}\la{test} $\Psi\in{\mathcal T}$ if there exist $N$, $w_1, \dots,
w_N \in C_0^{\infty}(\Rr^2)$, $\epsilon\ge 0$ and $\psi:\Rr^N\to \Rr$, smooth, such that
\[
\Psi(\theta) = \psi((J_{\epsilon}(\theta),w_1)_{\h},\dots, (J_{\epsilon}(\theta), w_N)_{\h})
\]
with $J_{\epsilon}$ a standard mollifier, i.e. convolution with $\epsilon^{-2}j(\fr{x}{\epsilon})$, $j\in C_0^{\infty}(\Rr^2)$, $j\ge 0$, $j(-x) = j(x)$, $\int_{\Rr^2} j(x)dx =1$, if $\epsilon>0$, and $J_{\epsilon} = {\mathbb I}$ if $\epsilon=0$.
\end{defi}
We note that the test functions are locally bounded and sequentially weakly continuous in $L^2(\Rr^2)$. We remind the elementary but important fact that weak continuity of real valued functions implies strong continuity, but in general continuity does not imply weak continuity. We identify $\Psi'(\theta)$ as an element of $L^2(\Rr^2)$ defined by 
\be
(\phi, \Psi'(\theta))_{\h} = 
\sum_{k=1}^N\left(\fr{\pa\psi}{\pa y_k}(y(\theta))\right)(J_{\epsilon}(\phi), w_k)_{\h}
\la{psiprime}
\ee
with
\be
y(\theta) = \left((J_{\epsilon}(\theta),w_1)_{\h},\dots, (J_{\epsilon}(\theta), w_N)_{\h}\right),
\la{ytheta}
\ee
that is
\be
\Psi'(\theta)(x) = \sum_{k=1}^N \left(\fr{\pa\psi}{\pa y_k}(y(\theta))\right) (J_{\epsilon}w_k)(x).
\la{psiprimepointwise}
\ee
We extend the definition (\ref{psiprime}) to more general $\phi$: this is the sense in which $(N^{(\nu)}(\theta),\Psi'(\theta))_{\h}$ is computed,
\be
(N^{(\nu)}(\theta),\Psi'(\theta))_{\h} = F_1(\theta) + \nu F_2(\theta) + F_3(\theta)
\la{nps}
\ee
with
\be
F_1(\theta) = \gamma(\theta, \D\Psi'(\theta))_{\h} - (f,\Psi'(\theta))_{\h},
\la{f1}
\ee
\be
F_2(\theta) = (\theta, (-\Delta)\Psi'(\theta))_{\h}
\la{f2}
\ee
and
\be
F_3(\theta) = -(\theta R^{\perp}\theta, \na\Psi'(\theta))_{\h}.
\la{f3}
\ee
Let us note that the Borel $\sigma$-algebra associated to the strong topology in $L^2(\Rr^2)$ is the same as the Borel $\sigma$ algebra associated to the weak topology because any open ball is a countable union of closed balls, which are convex, hence  weakly closed. The function $\theta\mapsto \|\theta\|_{H^1(\Rr^2)}^2$ is a Borel measurable function in $L^2(\Rr^2)$ because it is everywhere the limit of a sequence of continuous functions $\theta\mapsto \|J_{\epsilon}\theta\|_{H^{1}(\Rr^2)}^2$. The same of course applies to $\|\theta\|_{H^{\fr{1}{2}}(\Rr^2)}^2$. Therefore conditions (a) and (c) in Definition \ref{sssn} make mathematical sense. Moreover, condition (c) implies that $\mu^{(\nu)}$ is supported in the ball
\be
\|\theta\|_{H^{\fr{1}{2}}(\Rr^2)} \le \fr{1}{\gamma}\|f\|_{\h}
\la{h12bound}
\ee
as it is easily seen by taking $E_1\ge \gamma^{-1}\|f\|_{\h}$.
The integrand in condition (b) is locally bounded and weakly continuous:
\begin{lemma} \la{weakc} For any fixed $\Psi\in {\mathcal T}$ the maps
\[
\theta \mapsto F_i(\theta)
\]
$i=1,2,3$ are locally bounded in $L^2(\Rr^2)$ and weakly continuous in $L^2(\Rr^2)$ on bounded sets of $L^2(\Rr^2)\cap L^p(\Rr^2)$, $1\le p<2$. In particular
\[
\theta \mapsto (N^{(\nu)}(\theta),\Psi'(\theta))_{\h}
\]
is locally bounded and weakly continuous on bounded sets of  $L^2(\Rr^2)\cap L^p(\Rr^2)$, $1\le p<2$.
\end{lemma}
After inspection of the definitions (\ref{f1}, \ref{f2}, \ref{f3}) it is clear that the only nontrivial statement is about $F_3$. It is only for $F_3$ also that we need the bound in $L^p(\Rr^2)$ for $p<2$. If $\theta_n$ is weakly converging to $\theta$ in $L^2(\Rr^2)$, then the vectors $y(\theta_n)$ defined in (\ref{ytheta}) converge to the vector $y(\theta)$ and because the sequence $\theta_n$ is bounded in $L^2(\Rr^2)$, it follows that $y(\theta_n)$ belong to a fixed compact set in $\Rr^N$. The functions $\fr{\pa\psi}{\pa y_k}$ are continuous, so all we need to check is the convergence
\[
(\theta_nR^{\perp}\theta_n, \na J_{\epsilon}w_k)_{\h}\to (\theta R^{\perp}\theta, \na J_{\epsilon}w_k)_{\h}.
\]
In order to do this we make use of the assumed bound
\[
\sup_{n}\|\theta_n\|_{L^p(\Rr^2)}\le A_p.
\]
We note first that $\theta$, the weak limit in $L^2(\Rr^2)$, also obeys 
\[
\|\theta\|_{L^p(\Rr^2)}\le A_p.
\]
This follows by essentially restricting  $\theta^{p-1}$ (or ${\rm{sign}}\, \theta$ if $p=1$) on large compacts, integrating against $\theta_n$, passing to the limit in $n$ and then letting the compacts grow to the whole space. The weak continuity of the nonlinearity in SQG is proved using the representation ~\cite{Resnick95}
\be
(\theta R^{\perp}\theta, \na\phi)_{\h} = \fr{1}{2}\int_{\Rr^2}(\Lambda^{-1}\theta)(x)\left [\Lambda, \na\phi\right ](R^{\perp}(\theta)(x)dx 
\la{commutid}
\ee
for smooth compactly supported $\phi$,  with $[\Lambda, a]b = \Lambda(ab) -a\Lambda b$, the commutator of the operators $\Lambda$ and of multiplication by $a$. In our case $\phi= J_{\epsilon}w_k$.

In order to make sense of the terms in (\ref{commutid}) let us recall that 
the Riesz potential is given by ~\cite{Stein70} 
\[
\Lambda^{-1}\theta (x) = c\int_{\Rr^2}\fr{\theta(x-y)}{|y|}dy
\]
for an appropriate constant $c$. If $\theta \in L^p(\Rr^2)\cap L^{2}(\Rr^2)$ then
$\Lambda^{-1}\theta \in L^2(\Rr^2) + L^{\infty}(\Rr^2)$. Indeed, 
\[
\sup_{x}\left |\int_{|y|\ge 1}\fr{1}{|y|}\theta(x-y)dy\right |\le C\|\theta\|_{L^p(\Rr^2)}
\]
because $|y|^{-1}\in L^q(|y|\ge 1)$, $q>2$, $q^{-1}+p^{-1}=1$,
and
\[
\left \| \int_{|y|\le 1}\fr{1}{|y|}\theta(x-y)dz\right \|_{L^2(dx)}\le C\|\theta\|_{L^2(\Rr^2)}\]
as it is easily seen by duality or by Fourier transform.
Then, we note that
\[
\Lambda (ab)(x)-a(x)\Lambda b(x) = cP.V. \int_{\Rr^2}b(y)\fr{a(x)-a(y)}{|x-y|^3}dy
\]
and therefore, if $a$ is compactly supported in a ball of radius $L$ and if $|x| \ge 2L$, then, pointwise
\[
|\Lambda(ab)(x)-a(x)\Lambda b(x)| \le C|x|^{-3}\|a\|_{\h}\|b\|_{\h}.
\]
Thus, if $L$ is the radius of a ball in $\Rr^2$ containing the support of $\phi = J_{\epsilon}w_k$ and we denote
\[
C_{\phi}(\theta)(x) = \left [\Lambda, \na\phi\right]\cdot R^{\perp}\theta
\]
we have, for $\rho \ge 2L$,
\[
\left |\int_{|x|\ge \rho}(\Lambda^{-1}\theta )(x) C_{\phi}(\theta)(x)dx\right|
\le C\rho^{-1}\|\theta\|_{L^p(\Rr^2)}\|\theta\|_{L^2(\Rr^2)} + C\rho^{-2}\|\theta\|^2_{\h}
\]
We have thus
\[
\left |\int_{|x|\ge \rho}(\Lambda^{-1}\theta_n) (x) C_{\phi}(\theta_n)(x)dx\right|    \le \epsilon
\]
uniformly for all $n$ and also for $\theta$, provided $\rho$ is large enough 
so that $\rho\ge 2L$ and
\[
C\rho^{-1}A_pA_2 + C\rho^{-2}A_2^2 \le \epsilon
\]
where $A_2$ is the bound on $\|\theta_n\|_{\h}$.  It is well-known and easy to show that
\[
\theta \mapsto C_{\phi}(\theta) =\left[\Lambda, \na\phi\right ]R^{\perp}{\theta}
\]
is a bounded linear operator in $L^2(\Rr^2)$ for fixed $\phi\in C_0^{\infty}(\Rr^2)$. Thus,
\[
\left |\int_{|x|\le \rho}(\Lambda^{-1}\theta)(x)C_{\phi}(\theta)(x)dx\right| \le C\|\Lambda^{-1}\theta\|_{L^{2}(B(0,\rho))}\|\theta\|_{\h}
\]
The proof of the identity (\ref{commutid})  is best explained by denoting $\psi = \Lambda^{-1}\theta$ and $u = R^{\perp}\theta = \na^{\perp}\psi$. Then we have
\[
\ba
(\theta R^{\perp}\theta, \na\phi)_{\h} = \int_{\Rr^2}(\Lambda\psi)(u\cdot\na\phi)dx \\ = \int_{\Rr^2}\psi \left(\left [\Lambda, \na\phi\right ]\cdot u\right) dx 
+ \int_{\Rr^2}\psi \na\phi\cdot \Lambda u dx\\ = \int_{\Rr^2}\psi \left(\left [\Lambda, \na\phi\right ]\cdot u\right) dx
+ \int_{\Rr^2}\psi \na\phi\cdot \na^{\perp}\theta dx \\
= \int_{\Rr^2}\psi \left(\left [\Lambda, \na\phi\right ]\cdot u\right) dx -(\theta R^{\perp}\theta, \na\phi)_{\h}.
\ea
\]
In order to obtain the last term we integrated by parts and used $\na^{\perp}\cdot\na\phi = 0$. Now the weak continuity follows by writing
\[
\ba
\int_{\Rr^2}(\Lambda^{-1}\theta_n)(x) C_{\phi}(\theta_n)(x)dx - \int_{\Rr^2}(\Lambda^{-1}\theta)(x)C_{\phi}\theta(x)dx = \\
\int_{B(0,\rho)}(\Lambda^{-1}(\theta_n-\theta))(x) C_{\phi}(\theta_n)(x)dx +
\int_{B(0,\rho)}(\Lambda^{-1}\theta)(x) C_{\phi}(\theta_n -\theta)(x)dx \\
+ \int_{|x|\ge\rho}(\Lambda^{-1}\theta_n)(x) C_{\phi}(\theta_n)(x)dx
- \int_{|x|\ge \rho}(\Lambda^{-1}\theta)(x)C_{\phi}\theta(x)dx. 
\ea
\]
We pick $\epsilon>0$ and fix it. We choose $\rho>0$ large enough so that the last two terms are less than $\epsilon$ each. We fix $\rho$. The function
$(\Lambda^{-1}\theta) \chi_{B(0,\rho)}$ is a fixed function in $L^2(\Rr^2)$ (here $\chi_{B(0,\rho)}$ is the characteristic function), and, because $C_{\phi}$ is a bounded linear operator in $L^2(\Rr^2)$ the sequence $C_{\phi}(\theta_n-\theta)$ converges weakly to zero in $L^2(\Rr^2)$. Thus, letting $n\to \infty$ the ante-penultimate term converges to zero. Finally, for the first term
\[
\left |\int_{B(0,\rho)}(\Lambda^{-1}(\theta_n-\theta))(x) C_{\phi}(\theta_n)(x)dx \right| \le CA_2\|\Lambda^{-1}(\theta_n-\theta)\|_{L^2(B(0,\rho))}
\]
and this converges to zero because $\theta_n-\theta$ converges weakly to zero in $L^2(\Rr^2)$ and is bounded in $L^p(\Rr^2)$, $p<2$. Indeed, by the previous considerations about $\Lambda^{-1}$, for $\chi\in C_0^{\infty}(\Rr^2)$ we have that $\chi\Lambda^{-1}(\theta_n-\theta)$ is bounded in $H^1(\Rr^2)$ and converges weakly to 0 in $L^2(\Rr^2)$. Thus $\Lambda^{-1}(\theta_n-\theta)$ converges strongly to zero in $L^2(B(0,\rho))$ implying that the first term converges to zero. We conclude that the limit difference is in absolute value less than $2\epsilon$ and $\epsilon$ is arbitrary.

\begin{defi}\la{sssqg} A stationary statistical solution of the
forced critical SQG equation
\be
\pa_t \theta + R^{\perp}\theta\cdot\na\theta + \gamma \D\theta = f
\la{fsqg}
\ee
is a Borel probability measure $\mu$ on $L^2(\Rr^2)$ such that
\be 
\int_{\h}\|\theta\|^2_{H^{\fr{1}{2}}(\Rr^2)}d\mu(\theta)<\infty,
\la{h1finite}
\ee
and the equation
\be
\int_{\h} (N(\theta), \Psi'(\theta))_{\h}d\mu(\theta) = 0
\la{stateq}
\ee
holds for all cylindrical test functions $\Psi\in{\mathcal T}$, where 
\be
N(\theta) = R^{\perp}\theta\cdot\na\theta + \gamma\D\theta -f.
\la{ntheta}
\ee
We say that the stationary statistical solution satisfies the energy dissipation balance if
\be
\int_{\h}\left\{\gamma\|\theta\|^2_{H^{\fr{1}{2}}(\Rr^2)} - (f,\theta)_{\h}\right\}d\mu(\theta) = 0.
\la{edbal}
\ee 
\end{defi}
\begin{thm} \la{linv} Let  $\mu^{(\nu)}$ be a sequence of stationary statistical solutions of the viscous forced critical SQG equation (\ref{nfsqg}) with $f\in L^1(\Rr^2)\cap L^2(\Rr^2)$, with $\nu\to 0$. Assume that there exists $1\le p<2$ and 
$A_p$ such that the supports of the measures $\mu^{(\nu)}$ are included in
\be
B_p = \left\{\theta\in L^p(\Rr^2)\cap L^2(\Rr^2)\left |\right. \; \|\theta\|_{L^p(\Rr^2)}\le A_p\right\}.
\la{apball}
\ee
Then there exists a subsequence, denoted also $\mu^{(\nu)}$ and a stationary statistical solution $\mu$ of the forced critical SQG equation (\ref{fsqg}) such that
\be
\lim_{\nu\to 0}\int_{\h}\Phi(\theta)d\mu^{(\nu)}(\theta) = \int_{\h}\Phi(\theta)d\mu(\theta)
\la{limm}
\ee
holds for all weakly continuous, locally bounded real valued functions $\Phi$.
\end{thm}
As we saw in (\ref{h12bound}) the support of the measures $\mu^{(\nu)}$ is included in
\be
B = \left \{\theta\in H^{\fr{1}{2}}(\Rr^2)\left | \right.\; \|\theta\|_{H^{\fr{1}{2}}(\Rr^2)}\le \fr{\|f\|_{\h}}{\gamma}\right\}.
\la{bh12}
\ee
The set  $A_p =B\cap B_p$ is weakly closed in $L^2(\Rr^2)$ and it is a separable metrizable compact space with the weak $L^2(\Rr^2)$ topology. By Prokhorov's theorem the sequence $\mu^{(\nu)}$ is tight and therefore has a weakly convergent subsequence. The limit $\mu$  is a Borel probability on $A_p$. The extension of $\mu$ to $L^2(\Rr^2)$, denoted again by $\mu$ and given by $\mu(X) = \mu(X\cap A_p)$ is a Borel measure because $A_p$ is weakly closed. The measure $\mu$ satisfies (\ref{h1finite}) because it is supported in $B$. The equation (\ref{stateq}) is satisfied because we may pass to the limit in (2) of Definition \ref{sssn} in view of Lemma \ref{weakc}. 

\section{Inviscid Limit and Energy Dissipation Balance}
In this section we prove

\begin{thm}\la{limbal} Let $\mu^{(\nu)}$ be a sequence $\nu\to 0$ of stationary statistical solutions of the  forced viscous SQG equation (\ref{nfsqg}) supported in
\be
A = \left\{\theta\left |\right. \; \|\theta\|_{L^p(\Rr^2)}\le A_p,\; \|\theta\|_{L^{\infty}(\Rr^2)}\le A_{\infty},\;
 \|\theta\|_{H^{\fr{1}{2}}(\Rr^2)}\le \fr{\|f\|_{\h}}{\gamma}\right\}.
\la{abound}
\ee
Let $\mu$ be any weak limit of $\mu^{(\nu)}$ in $L^2(\Rr^2)$. Then $\mu$ is a stationary statistical solution of the forced critical SQG equation (\ref{fsqg}) that satisfies the energy dissipation balance (\ref{edbal}).
\end{thm} 
In fact, by Theorem \ref{linv}, we know that any weak limit is a stationary statistical solution of the forced critical SQG equation. We check that it is supported on $A$. The  set $A$ is weakly closed in $L^2(\Rr^2)$, and because its complement $U$ is weakly open and
\[
\mu(U)\le \lim\inf_{\nu\to 0}\mu^{(\nu)}(U)=0
\]
it follows that $\mu$ is supported in the set $A$. The rest of the proof is
done by showing that (\ref{stateq}) and the fact that $\mu$ is supported in $A$ imply (\ref{edbal}).

We take a sequence $w_j\in C_0^{\infty}(\Rr^2)$ that is an orthonormal basis of $L^2(\Rr^2)$. We fix $\epsilon>0$ and consider the sequence of test functions
\[
\Psi_m(\theta) = \fr{1}{2}\sum_{k=1}^m (J_{\epsilon}(\theta), w_j)^2_{\h}
\]
i.e. we take $\psi(y) = \fr{1}{2}\sum_{k=1}^my_k^2$ in Definition \ref{test}.
We note that
\[
(N(\theta),\Psi'_m(\theta))_{\h} = \sum_{k=1}^m(J_{\epsilon}(\theta),w_j)_{\h}(J_{\epsilon}(N(\theta)),w_j)_{\h}.
\]
This is uniformly bounded in $m$ because
\[
\left |(N(\theta),\Psi'_m(\theta))_{\h}\right| \le \|J_{\epsilon}\theta\|_{\h}\|J_{\epsilon}N(\theta)\|_{\h}.
\]
On $A$, the right hand side is bounded uniformly in $\theta$
\[
\left |(N(\theta),\Psi'_m(\theta))_{\h}\right| \le A_2((2+\epsilon^{-\fr{1}{2}})\|f\|_{\h} + \epsilon^{-1}A_2A_{\infty})
\]
and, by Parseval, it converges to $(J_{\epsilon}(\theta),J_{\epsilon}(N(\theta)))_{\h}$ pointwise, as $m\to\infty$. So, we deduce from (\ref{stateq}) and the Lebesgue dominated convergence theorem that
\[
\int_{\h}(J_{\epsilon}(\theta), J_{\epsilon}(N(\theta)))_{\h}d\mu(\theta) = 0
\]
for any  $\epsilon>0$. This can be written as
\be
I_{\epsilon} + K_{\epsilon} = 0
\la{iepluske}
\ee
where the two terms are
\be
I_{\epsilon} = \int_{\h}(J_{\epsilon}(\theta), J_{\epsilon}(\gamma\D\theta -f))_{\h}d\mu(\theta)
\la{ie}
\ee
and
\be
K_{\epsilon} = \int_{\h}(J_{\epsilon}(\theta), J_{\epsilon}((R^{\perp}\theta)\cdot\na\theta))_{\h}d\mu(\theta).
\la{ke}
\ee
Now
\[ 
(J_{\epsilon}(\theta), J_{\epsilon}(u)\cdot\na J_{\epsilon}(\theta))_{\h} =0
\]
in view of the incompressibility of $u =R^{\perp}{\theta}$, so
\be
K_{\epsilon} = \int_{\h}(J_{\epsilon}(\theta), \na\cdot\rho_{\epsilon}(u,\theta))_{\h}d\mu(\theta).
\la{kenew}
\ee
where
\be
\rho_{\epsilon}(u,\theta) = J_{\epsilon}(u\theta) -(J_{\epsilon}(u))(J_{\epsilon}(\theta)). 
\la{rhodef}
\ee
We show that $\lim_{\epsilon\to 0}K_{\epsilon} =0$. We use the identity ~\cite{ConstantinETiti94}
\[
\rho_{\epsilon}(u,\theta) = r_{\epsilon}(u,\theta) - (u-J_{\epsilon}(u))(\theta-J_{\epsilon}(\theta))
\]
with
\[
r_{\epsilon}(u,\theta)(x) = \int_{\Rr^2}j(z)(\delta_{\epsilon z}(u)(x))(\delta_{\epsilon z}(\theta)(x))dz,
\]
\[
\delta_z(u)(x) = u(x-\epsilon z)-u(x)
\]
and
\[
\delta_{\epsilon z}(\theta)(x) = \theta(x-\epsilon z)-\theta(x).
\] 
Clearly also
\[
(J_{\epsilon}u -u)(J_{\epsilon}\theta-\theta) = \int_{\Rr^4}j(z_1)j(z_2)(\delta_{\epsilon z_1}u) (\delta_{\epsilon z_2}\theta) dz_1dz_2.
\]
The inequality 
\be
\|\delta_{\epsilon z} \theta \|_{L^4(\Rr^2)}^2 \le C(\epsilon |z|)^{\fr{1}{2}}\|\theta\|_{L^{\infty}(\Rr^2)}\|\theta\|_{H^{\fr{1}{2}}(\Rr^2)}.
\la{l4dtheta}
\ee
and its consequence (because of the boundeness of $R^{\perp}$ in $L^4(\Rr^2)$)
\be
\|\delta_{\epsilon z} R^{\perp}\theta \|_{L^4(\Rr^2)}^2 \le C(\epsilon |z|)^{\fr{1}{2}}\|\theta\|_{L^{\infty}(\Rr^2)}\|\theta\|_{H^{\fr{1}{2}}(\Rr^2)}
\la{l4du}
\ee
follow from the elementary inequality
\be
\|\delta_{\epsilon z}\theta\|_{L^2(\Rr^2)}\le C(\epsilon |z|)^{\fr{1}{2}}\|\theta\|_{H^{\fr{1}{2}}(\Rr^2)}
\la{l2dtheta}
\ee
which is proved by Fourier transform. Consequently,
\be
\|\rho_{\epsilon}(R^{\perp}\theta_1,\theta_2)\|_{\h}\le C \epsilon^{\fr{1}{2}}\|\theta_1\|_{L^{\infty}(\Rr^2)}^{\fr{1}{2}}\|\theta_2\|_{L^{\infty}(\Rr^2)}^{\fr{1}{2}}\|\theta_1\|_{H^{\fr{1}{2}}(\Rr^2)}^{\fr{1}{2}}\|\theta_2\|_{H^{\fr{1}{2}}(\Rr^2)}^{\fr{1}{2}}.
\la{rhobound}
\ee
The integrand in $K_{\epsilon}$ is bounded 
\be
\left |(\na J_{\epsilon}(\theta), \rho_{\epsilon}(u,\theta))_{\h}\right| \le C
\|\theta\|_{L^{\infty}(\Rr^2)}\|\theta\|_{H^{\fr{1}{2}}(\Rr^2)}^2
\la{kep}
\ee
and converges to $0$ as $\epsilon\to 0$, for fixed $\theta\in H^{\fr{1}{2}}(\Rr^2)\cap L^{\infty}(\Rr^2)$. Indeed, the trilinear map
\[
(\theta_1,\theta_2,\theta_3) \mapsto (\na J_{\epsilon}(\theta_1), \rho_{\epsilon}(R^{\perp}\theta_1,\theta_2))_{\h}
\]
is separately continuous from $H^{\fr{1}{2}}(\Rr^2)$ to $\Rr$ uniformly on $A$ and uniformly in $\epsilon$. This can be seen from the expression
\[
 (\na J_{\epsilon}(\theta_3), \rho_{\epsilon}(R^{\perp}\theta_1, \theta_2))_{\h}=
-\fr{1}{\epsilon}\int\na_zj(z)(\delta_{\epsilon z}\theta_3, \rho_{\epsilon}(R^{\perp}\theta_1, \theta_2))_{\h}dz
\]
and the bound obtained from (\ref{l2dtheta}) and (\ref{rhobound})  
\be
\ba
\left | (\na J_{\epsilon}(\theta_3), \rho_{\epsilon}(R^{\perp}\theta_1, \theta_2))_{\h}\right |\\
\le C\|\theta_3\|_{H^{\fr{1}{2}}(\Rr^2)}\|\theta_1\|_{L^{\infty}(\Rr^2)}^{\fr{1}{2}}\|\theta_2\|_{L^{\infty}(\Rr^2)}^{\fr{1}{2}}\|\theta_1\|_{H^{\fr{1}{2}}(\Rr^2)}^{\fr{1}{2}}\|\theta_2\|_{H^{\fr{1}{2}}(\Rr^2)}^{\fr{1}{2}}.
\ea
\la{trib}
\ee
This explains (\ref{kep}) and also shows the pointwise vanishing of the integrand in $K_{\epsilon}$ as $\epsilon\to 0$: the integrand in (\ref{kenew}) obviosly tends to zero for smooth $\theta$, and $\theta$ in $A$ can be approximated in the norm of $H^{\fr{1}{2}}(\Rr^2)$ by smooth functions.
Therefore, from the Lebegue dominated convergence theorem
\[
\lim_{\epsilon\to 0}K_{\epsilon} =0.
\]
It remains to consider the limit of $I_{\epsilon}$, but this is quite straightforward,
\[
\lim_{\epsilon\to 0}I_{\epsilon} = \int_{\h}(\gamma\|\theta\|_{H^{\fr{1}{2}}(\Rr^2)}^2 - (\theta, f)_{\h})d\mu(\theta) 
\]
and thus the proof is complete.

\section{Long Time Averages}
In this section we consider long time averages of solutions and the stationary statistical solutions they generate. We start by recalling the concept of generalized (Banach) limit:
\begin{defi}\la{gl} A generalized limit (Banach limit) is a bounded linear functional
\[
\bl: BC([0,\infty))\to \Rr
\]
such that 

1. $\bl(g)\ge 0, \quad \forall g\in BC([0,\infty)), \; g\ge 0.$

2. $\bl(g) = \lim_{t\to\infty}g(t)$ whenever the usual limit exists.
\end{defi}
The space $BC([0,\infty))$ is the Banach space of all bounded continuous real valued functions defined on $[0, \infty)$ endowed with the $\sup$ norm. It can be shown that the generalized limit always saitisfies
\[
\lim\inf_{t\to\infty}g(t)\le \bl (g)\le \lim\sup_{t\to\infty} g
\]
for all $g\in BC([0,\infty))$. Moreover, given a fixed $g_0\in BC([0,\infty))$ and a sequence $t_j\to\infty$ for which $\lim_{j\to\infty}g_0(t_j) = l$ exists,
a generalized limit can be found which satisfies $\bl(g_0) = l$. This implies that one can choose a generalized limit that obeys $\bl(g_0) = \lim\sup_{t\to\infty}g_0(t)$. 
With this language at our disposal, we now state the result about long time averages of viscous forced critical SQG.
\begin{thm} \la{bogo} Let $f\in L^1(\Rr^2)\cap L^{\infty}(\Rr^2)$ and $\theta_0\in  L^1(\Rr^2)\cap L^{\infty}(\Rr^2)$. Let $\bl$ be a Banach limit. Then the map
\be
\Phi \mapsto \bl\fr{1}{t}\int_0^t\Phi(S^{(\nu)}(s,\theta_0))ds
\la{limap}
\ee
for $\Phi \in C(L^2(\Rr^2))$  (strongly continuous, real valued) defines a a stationary statistical solution $\mu^{(\nu)}$ of the viscous forced SQG equation (\ref{nfsqg}):
\be
\int_{\h}\Phi(\theta)d\mu^{(\nu)}(\theta) = \bl\fr{1}{t}\int_0^t\Phi(S^{(\nu)}(s,\theta_0))ds.
\la{munulim}
\ee
The measure is supported in the set
\be
A = \left\{\theta\left| \right. \|\theta\|_{H^{\fr{1}{2}}(\Rr^2)} \le \fr{\|f\|_{\h}}{\gamma}, \; \|\theta\|_{L^p(\Rr^2)}\le A_p,\; 1\le p\le \infty\right\}
\la{aset}
\ee
with 
\[
A_p =  \|\theta_0\|_{L^p(\Rr^2)} + \fr{\|f\|_{L^p(\Rr^2)}}{\gamma},\; 1\le p\le \infty.
\]
The inequality
\be
\int_{\h}\left [\nu\|\na\theta\|^2_{\h} + \gamma\|\theta\|_{H^{\fr{1}{2}}(\Rr^2)}^2 - (f, \theta)_{\h}\right ]d\mu^{(\nu)}(\theta) \le 0
\la{enineq}
\ee
holds.
\end{thm}
The positive semiorbit 
\[
O_{+}(\theta_0) = \{\theta\left |\right. \; \exists s\ge 0, \; \theta = S^{(\nu)}(s, \theta_0)\}
\]
is relatively compact in $L^2(\Rr^2)$ because it is bounded in $H^{1}(\Rr^2)$ and uniformly integrable by Proposition \ref{nsprop}, (\ref{notr}). For any 
$\Phi\in C\left({\overline{O_{+}(\theta_0)}}\right )$ the function $s\mapsto \fr{1}{t}\int_0^t\Phi(S^{(\nu)}(s,\theta_0))ds$ is a bounded continuous function on $[0,\infty)$ so we may apply $\bl$ to it. (Of course, $C(L^2(\Rr^2))\subset C\left({\overline{O_{+}(t_0, \theta_0)}}\right)$.) The map
\[
\Phi\mapsto \bl\fr{1}{t}\int_0^t\Phi(S^{(\nu)}(s,\theta_0))ds
\]
is a positive functional on $C\left(\overline{O_{+}(\theta_0)}\right )$. Because of the Riesz representation theorem on compact spaces, it follows that there exists a Borel measure representing it, i.e. (\ref{munulim}) holds. The measure is supported on $\overline{O_{+}(\theta_0)}$. We take a test function $\Psi\in{\mathcal T}$. Then
\[
\int_{\h}(N^{(\nu)}(\theta), \Psi'(\theta))_{\h}d\mu^{(\nu)}(\theta) =
\bl \fr{1}{t}\int_0^t\fr{d}{ds}\Psi(S^{(\nu)}(s,\theta_0))ds 
\]
holds and the right hand side vanishes, verifying  (b) of Definition (\ref{sssn}). Because of (\ref{lpb}) the semiorbit is included in the set
\[
\{\theta\left |\right. \; \|\theta \|_{L^p(\Rr^2)} \le A_p, \;\; 1\le p\le \infty\}.
\]
The fact that the support of $\mu^{(\nu)}$ is included in $A$ follows as shown before  from property (c) of Definition \ref{sssn}. 
In order to check (a), (c) of Definition \ref{sssn} we would like to take long time averages in the energy balance (\ref{bal}). In order to do so, we first mollify the equation.  This is due to the fact that $\|\na\theta\|_{L^2(\Rr^2)}^2$ is not continuous in $L^2(\Rr^2)$. We put
\[
\theta_{\epsilon} (x,t) = J_{\epsilon}(S^{(\nu)}(t,\theta_0))(x),\quad u_{\epsilon}(x,t) = J_{\epsilon}R^{\perp}(S^{(\nu)}(t,\theta_0)),
\]
and applying $J_{\epsilon}$ to (\ref{nfsqg}), multiplying by $\theta_{\epsilon}$ and integrating we deduce
\[
\ba
\fr{1}{t}\int_0^t \left [\gamma \|\theta_{\epsilon}(s)\|^2_{H^{\fr{1}{2}}(\Rr^2)} - (J_{\epsilon}f, \theta_{\epsilon}(s))_{\h} + \nu\|\na\theta_{\epsilon}(s)\|^2_{\h}\right]ds \\=
\fr{1}{2t}\left[\|\theta_{\epsilon}(0)\|_{\h}^2 - \|\theta_{\epsilon}(t)\|_{\h}^2\right] +
\fr{1}{t}\int_0^t(\rho (u_{\epsilon}(s), \theta_{\epsilon}(s)), \na\theta_{\epsilon}(s))_{\h}ds.
\ea
\]
We obtain
\be
\ba
\int\left [\gamma\|J_{\epsilon}\theta\|_{H^{\fr{1}{2}}(\Rr^2)}^2 - (J_{\epsilon}f, J_{\epsilon}\theta)_{\h} + \nu\|\na J_{\epsilon}\theta\|_{\h}^2\right ]d\mu^{(\nu)}(\theta)\\
= \bl\fr{1}{t}\int_0^t(\rho (u_{\epsilon}(s), \theta_{\epsilon}(s)), \na\theta_{\epsilon}(s))_{\h}ds.
\ea
\la{ltineq}
\ee
Because of (\ref{bal}) and (\ref{lpb})
\be
\ba
\lim\sup_{t\to\infty}\fr{1}{t}\int_0^t\left [\gamma\|S^{(\nu)}(s, \theta_0)\|_{H^{\fr{1}{2}}(\Rr^2)}^2 + \nu\|\na S^{(\nu)}(s,\theta_0)\|^2_{L^2(\Rr^2)}ds\right ] \\
\le \fr{1}{\gamma}\|f\|_{\h}^2
\ea
\la{nltb}
\ee
and because $J_{\epsilon}$ does not increase $L^2$ norms, we deduce from (\ref{nltb}) that
\[
\sup_{0<\epsilon}\int_{\h}\left[\gamma\|J_{\epsilon}\theta\|_{H^{\fr{1}{2}}(\Rr^2)}^2 + \nu\|\na J_{\epsilon}\theta\|_{L^2(\Rr^2)}^2 \right]d\mu^{(\nu)}(\theta) \le  \fr{1}{\gamma}\|f\|_{\h}^2.
\]
The functions $\|\theta\|_{H^{\fr{1}{2}}(\Rr^2)}^2$ and   $\|\na \theta\|_{L^2(\Rr^2)}^2$ are Borel measurable and so, by Fatou, we obtain (a) of Definition \ref{sssn}.
Using the $H^1\cap L^{\infty}$ information we have
\[
\|\rho_{\epsilon}(R^{\perp}\theta, \theta)\|_{L^2(\Rr^2)} \le C\sqrt{\epsilon}\|\theta\|_{L^{\infty}(\Rr^2)}\|\na\theta\|_{L^2(\Rr^2)}
\]
and thus
\[
\ba
\bl\fr{1}{t}\int_0^t\rho (u_{\epsilon}(s)\theta_{\epsilon}(s)), \na\theta_{\epsilon}(s))_{\h}ds\\
\le C\epsilon\left [\|\theta_0\|_{L^{\infty}(\Rr^2)} + \fr{1}{\gamma}\|f\|_{L^{\infty}(\Rr^2)}\right ]\fr{1}{\nu\gamma}\|f\|_{\h}^2.
\ea
\]
This implies that the right hand side of (\ref{ltineq}) converges to zero as $\epsilon\to 0$. This proves (\ref{enineq}) by Fatou. In order to prove (c) of Definition \ref{sssn} we take $\chi'(y)$, a smooth, nonnegative, compactly supported function defined for $y\ge 0$. Then $\chi(y) = \int_0^y\chi'(e)de$ is bounded on $\Rr_+$ and 
\[
\fr{d}{dt}\chi(\|\theta_{\epsilon}(t)\|_{\h}^2) = 
\chi'(\|\theta_{\epsilon}(t)\|_{\h}^2)\fr{d}{dt}\|\theta_{\epsilon}(t)\|^2_{\h}.
\]
We proceed as above and obtain
\[
\ba
\int_{\h}\chi'(\|\theta\|^2_{\h})\left\{\nu\|\na\theta\|_{\h}^2 + \gamma\|\theta\|_{H^{\fr{1}{2}}(\Rr^2)}^2 - (f,\theta)_{\h}\right\}d\mu^{(\nu)}(\theta)\\
\le 0.
\ea
\]
We let $\chi'(y)$ converge pointwise to the characteristic function of the interval $[E_1^2, E_2^2]$ with $0\le \chi'(y)\le 2$. This proves (c) of Definition \ref{sssn} and concludes the proof of this theorem.
\section{Conclusion}
\begin{thm} Let $\theta_0, f \in L^1(\Rr^2)\cap L^{\infty}(\Rr^2)$. Then
\[
\lim_{\nu\to 0}\nu\left(\limsup_{t\to\infty}\fr{1}{t}\int_0^t\|\na S^{(\nu)}(s,\theta_0)\|^2_{\h}ds\right) = 0.
\]
\end{thm}
We argue by contradiction. If the conclusion would be false, then there would 
exist $\delta>0$, a sequence $\nu_k\to 0$, and, for each $\nu_k$, a sequence of times $t_j\to \infty$ such that
\[
\fr{\nu_k}{t_j}\int_0^{t_j}\|\na S^{(\nu_k)}(s,\theta_0)\|^2_{\h}ds \ge \delta
\]
holds for all $t_j$. Because of (\ref{bal})
\[
\ba
\delta\le \fr{\nu_k}{t_j}\int_0^{t_j}\|\na S^{(\nu_k)}(s,\theta_0)\|^2_{\h}ds =\\
\fr{1}{t_j}\int_0^{t_j}\left[-\gamma\|S^{(\nu_k)}(s,\theta_0)\|^2_{H^{\fr{1}{2 }}(\Rr^2)} + (f, S^{(\nu_k)}(s,\theta_0))_{\h}\right]ds \\
+ \fr{1}{2t_j}\left [\|\theta_0\|^2_{\h} - \|S^{(\nu_k)}(t,\theta_0)\|^2_{\h}\right]
\ea
\]
It follows that
\be
\limsup_{t\to\infty}\fr{1}{t}\int_0^t\left[-\gamma\|S^{(\nu_k)}(s,\theta_0)\|^2_{H^{\fr{1}{2 }}(\Rr^2)} + (f, S^{(\nu_k)}(s,\theta_0))_{\h}\right]ds \ge   \delta.
\la{limsupdelta}
\ee
By Theorem \ref{bogo} there exists a stationary statistical solution of the forced viscous SQG equation, $\mu^{(\nu_k)}$ supported in $A$ such that
\be
\int_{\h}\left\{-\gamma\|\theta\|^2_{H^{\fr{1}{2}}(\Rr^2)} + (f,\theta)_{\h}\right\}d\mu^{(\nu_k)}(\theta)\ge \delta>0.
\la{mukdel}
\ee
Passing to a weakly convergent subsequence (denoted again $\mu^{(\nu_k)}$), we find using Theorem {\ref{linv}} and Theorem {\ref{limbal}} a stationary statistical solution $\mu$ of the forced critical SQG equation that satisfies the energy dissipation balance (\ref{edbal}). The function $\theta\mapsto (f,\theta)_{\h}$ is weakly continuous, so
\[
\lim_{k\to\infty}\int_{\h} (f,\theta)d\mu^{(\nu_k)}(\theta) = \int_{\h} (f,\theta)d\mu.
\]
On the other hand, by Fatou
\[
\int_{\h}\|\theta\|^2_{H^{\fr{1}{2}}}d\mu(\theta) \le \liminf_{k\to\infty}\int_{\h}\|\theta\|^2_{H^{\fr{1}{2}}}d\mu^{(\nu_k)}(\theta).
\]
Using (\ref{mukdel}) we obtain
\[
\int_{\h}\left[\gamma \|\theta\|^2_{H^{\fr{1}{2}}} - (f,\theta)_{\h}\right ]d\mu(\theta) \le -\delta<0
\]
contradicting (\ref{edbal}). This concludes the proof of the theorem.

The forced critical SQG equation is dissipative, and the main result here shows that additional viscous dissipation does not leave any anomalous remanent dissipation. The same result is true for spatially periodic boundary conditions, and for additional dissipation of the type $\nu(-\Delta)^{\alpha}$. The problem of absence of anomalous dissipation is open for the forced SQG equation without the $\Lambda$ term in $\D$, i.e. with friction that does not grow like $|k|$ for high wave-numbers $k$. 

The method of proof of ~\cite{ConstantinRamos07} and of this paper is quite general, and is applicable for a large class of equations where no uniform bound on the dissipation is readily available. The main ingredients necessary for the success of the method are: an energy dissipation balance for viscous solutions, relative compactness of viscous semiorbits, weak continuity of the nonlinearity, and enough bounds to control the nonlinear fluxes. The forced SQG equation with wave-number independent friction and the supercritical forced SQG equation have all the mentioned ingredients, except the last one, so what is missing is proving the energy dissipation balance for the long time averages of solutions of the inviscid equation.

{\bf Acknowledgments} The work of PC was supported in part by NSF grants DMS-1209394, DMS-1265132,  and DMS-1240743. The work of VV was supported in part by the NSF grant DMS-1211828.

\end{document}